\documentclass{amsart}

\usepackage{amsmath, amssymb,amsfonts,amsthm}

\theoremstyle{plain}
\newtheorem{thm}{Theorem}[section]
\newtheorem{propn}[thm]{Proposition}
\newtheorem{lemma}[thm]{Lemma}
\newtheorem{coroll}[thm]{Corollary}
\theoremstyle{definition}
\newtheorem{defn}[thm]{Definition}
\newtheorem{eg}[thm]{Example}
\theoremstyle{remark}
\newtheorem*{rem}{Remark}

\newcommand{\mbb}{\mathbb}
\newcommand{\mc}{\mathcal}

\newcommand{\lm}{\lambda}
\newcommand{\sig}{\sigma}
\newcommand{\tht}{\theta}
\newcommand{\Om}{\Omega}

\newenvironment{inline}[1]
{\vspace{0.5em} \noindent {\bf {#1}} }{\vspace{0.5em}}

\newcommand{\st}{\;:\;}
\newcommand{\defeq}{:=}
\newcommand{\dt}[1]{{\sf {#1}}}
\renewcommand{\emph}[1]{{\sl #1}\/}

\newcommand{\pcat}[1]{{\sf {#1}}} 

\newcommand{\id}[1][]{{\sf 1}_{#1}} 
\newcommand{\vctR}[1]{{\bf {#1}}} 
\newcommand{\blob}{\bullet}
\newcommand{\blank}{\underline{\quad}}
\newcommand{\lub}{\vee}

\newcommand{\iso}{\cong}

\newcommand{\tp}{\otimes}
\newcommand{\ptp}{\widehat{\otimes}} 
\newcommand{\ptpR}[1]{\underset{#1}{\ptp}} 
\newcommand{\lHom}[1]{{}_{#1}{\rm Hom}} 

\newcommand{\abs}[1]{\vert{#1}\vert}

\newcommand{\norm}[1]{\Vert{#1}\Vert}

\newcommand{\Cplx}{\mbb C}
\newcommand{\Nat}{\mbb N}
\newcommand{\Rat}{\mbb Q}
\newcommand{\Zahl}{\mbb Z}

\DeclareMathOperator{\Ker}{Ker}
\DeclareMathOperator{\Image}{Im}

\newcommand{\ackndiag}{Uses Paul Taylor's {\tt diagrams.sty} macros.}

\usepackage[PostScript=dvips]{diagrams} 
\newcommand{\LR}[2]{\pile{\lTo^{#1} \\ \rTo_{#2} }}

\newarrow{Dbl} ===={=>} 

\DeclareMathOperator{\Tor}{Tor}
\DeclareMathOperator{\Ext}{Ext}
\newcommand{\bdy}{{\sf d}}
\newcommand{\dif}{\delta}

\newcommand{\base}[3]{\mc{#1}^{#2}_{#3}}  
\newcommand{\Ho}[2]{\base{#1}{}{#2}}  
\newcommand{\Co}[2]{\base{#1}{#2}{}}  

\newcommand{\lp}[1]{\ell^{#1}}
\newcommand{\A}[1]{{\mc A}_{#1}} 
\newcommand{\fu}[1]{{#1}^{\#}} 

\newcommand{\SHo}[3]{\Ho{#1}{#2}({#3})}
\newcommand{\SCo}[3]{\Co{#1}{#2}({#3})}

\newcommand{\Ban}{\pcat{Ban}}
\newcommand{\BAlg}{\pcat{BAlg}}
\newcommand{\SGp}{\pcat{SGp}}
\newcommand{\SLat}{\pcat{SLat}}

\newcommand{\subst}[1]{\widehat{\vctR{#1}}}
\newcommand{\higheq}[1]{\boxed{ {#1} }}

\begin{document}
\title[Vanishing cohomology of Banach semilattice algebras]{Simplicial homology and Hochschild cohomology of Banach semilattice algebras}
\author{Yemon Choi}
\address{School of Mathematics and Statistics, University of Newcastle, Newcastle upon Tyne, NE1 7RU, England}
\email{y.choi.97@cantab.net}
\subjclass{Primary {\bf ???}; Secondary {\bf ???}}
\thanks{\ackndiag}
\begin{abstract}
The $\lp{1}$-convolution algebra of a semilattice is known to have trivial cohomology in degrees 1,2 and 3 whenever the coefficient bimodule is symmetric. We extend this result to all cohomology groups of degree $\geq 1$ with symmetric coefficients. Our techniques prove a stronger splitting result, namely that the splitting can be made natural with respect to the underlying semilattice.
\end{abstract}
\maketitle

\begin{section}{Introduction}
Let $S$ be a \emph{commutative} semigroup and $\A{S}=\lp{1}(S)$ its associated Banach convolution algebra. Gr{\o}nb{\ae}k proved in \cite{Gronbaek} that $\A{S}$ is amenable if and only if $S$ is a finite semilattice of abelian groups. Rather than asking about amenability, we can investigate when the Hochschild cohomology groups $\Co{H}{n}(\A{S},M)$ vanish for various $n$ and specific classes of $\A{S}$-bimodules $M$. Unfortunately, calculating $\Co{H}{2}(\A{S},M)$ has proved hard even when $M$ is a symmetric $\A{S}$-module.

Progress has been made in the case where $S$ is a \dt{semilattice}, that is, a commutative semigroup where every element is idempotent. If $S$ is a semilattice, then by an old result of Duncan and Namioka $\A{S}$ is amenable if and only if $S$ is finite (see \cite{Dunc-Nam}, Theorem 10). Nevertheless, Dales and Duncan observed in \cite{Da-Dunc} that $\Co{H}{1}(\A{S},M)=\Co{H}{2}(\A{S},M)=0$ for all semilattices $S$ and all symmetric $\A{S}$-bimodules $M$, and this has been extended to the third cohomology groups in \cite{GoPoWh}.

In this note we extend these results to all higher cohomology groups (the precise statement will be given in Section \ref{s:state_thm}). The actual calculations used to prove our main result generalise  familiar {\it ad hoc} arguments; our contribution is to put these arguments into a systematic framework, by using some basic language and ideas from category theory to ``formalise'' the process for solving the cohomology problems.

Let us review the structure of this paper. For readers unfamiliar with semigroups and their convolution algebras, we give the basic definitions and some examples. This is followed by a brief resum\'e on Banach algebras and their cohomology theory.

Our main result (Theorem \ref{t:result}) is stated in Section \ref{s:state_thm}; a special case is established in Section \ref{s:finite_free}, which will turn out to be the kernel for the proof of the general case. Section \ref{s:formalsub} contains the algebraic machinery which will be used to boost the results of Section \ref{s:finite_free} up to the general case, and Section \ref{s:mainthm} assembles the preceding work to give the main result. It is slightly unusual that the final proof is intrinsically inductive (on the degree of homology).
\end{section}

\begin{section}{Background and definitions}\label{s:background}
We refer to standard references such as \cite{Bons-Dunc} for the definitions of Banach algebras and bimodules over them. In particular, recall that if $A$ is a Banach algebra with identity $\id[A]$ and $M$ a Banach $A$-bimodule, then $M$ is said to be \dt{unit-linked} (\cite{Bons-Dunc}, \S9, Defn 11) if
\[ \id[A]\cdot x = x = x\cdot\id[A] \quad\quad\text{ for all $x \in M$}. \]

\begin{subsection}{Semigroup algebras and monoid algebras}
Let $S$ be a semigroup. We denote by $\A{S}$ the \dt{$\lp{1}$-convolution algebra of $S$} (for the full definition, see e.g. \cite{Bons-Dunc} \S1, Example 23). This construction is functorial: more precisely, $\A{\blob}$ defines a functor from the category $\SGp$ of semigroups and homomorphisms to the category $\BAlg_1$ of Banach algebras and contractive algebra homomorphisms.

Note that the class of algebras $\A{S}$ contains more familiar classes of Banach algebras as special cases:
\begin{itemize}
\item[--] if $S$ is a group then $\A{S}$ is just the usual $\lp{1}$-group algebra;
\item[--] if $S=\Zahl_+^k$ with the usual additive operations, then $\A{S}=\lp{1}(\Zahl_+^k)$ may be identified with the $\lp{1}$-completion of the polynomial ring $\Cplx[z_1,\ldots,z_k]$.
\end{itemize}
We also remark that our semigroups need not have identity elements: a semigroup with an identity element is called a \dt{monoid}.

This paper is concerned with a special class of commutative semigroups, basic in semigroup theory but perhaps less familiar in the context of functional analysis.
\begin{defn}
Let $S$ be a semigroup. An element $s$ of $S$ is \dt{idempotent} if $s^2=s$. A commutative semigroup $S$ in which every element is idempotent is called a \dt{semilattice}; if $S$ has an identity element then we shall say that $S$ is a \dt{unital semilattice}.
\end{defn}

\begin{eg}\label{eg:finfree}
Let $X$ be any set and consider the power set $2^X$. The binary map $2^X \times 2^X \to 2^X$ given by $(A,B) \mapsto A\cup B$ is associative, and one may therefore regard $2^X$ as a semigroup with product given by $\cup$. Clearly this semigroup is commutative, has a unit (namely the empty set $\emptyset$) and consists of idempotents, and thus $2^X$ is a unital semilattice -- in fact, it is the \dt{free unital semilattice} on the generating set $X$.
\end{eg}

\begin{eg}
Define $\lub: \Rat\times \Rat \to \Rat$ by $(m,n) \mapsto \max(m,n)$. Then $\lub$ defines an associative binary product, and the semigroup $(\Rat, \lub)$ turns out to be a non-unital semilattice. (More generally, any totally ordered set $T$ is a semilattice when equipped with the multiplication $(x,y)\mapsto x\vee y$.)
\end{eg}

\begin{rem}
Every commutative semigroup $S$ admits a decomposition $S=\coprod_{e \in L} S_e$ where $L$ is the set of idempotents in $S$ and $S_e\cdot S_f \subseteq S_{ef}$ (see, for example, \cite{How_fund-sgp}, Ch. 4, Exercise 3.). This suggests the possibility of relating homological properties of $\A{S}$ to those of each subalgebra $\A{S_e}$. One possible first step is to understand as much as possible about the case where each $S_e$ consists of a single element, i.e. when $S=L$ is a semilattice; the present paper is an attempt to make progress in this direction.
\end{rem}
\end{subsection}

\subsection*{Working with tensor powers}
In performing our main calculations we shall need to deal with projective tensor powers of $\A{S}$ (see \cite{Ryan_TP}, Ch. 2 for a full discussion of the \dt{projective tensor product} of Banach spaces). Throughout this paper the projective tensor product of two Banach spaces $E$ and $F$ is denoted by $E\ptp F$.

Potential difficulties with the projective tensor product can be sidestepped in what follows, because the underlying Banach space of $\A{S}$ is $\lp{1}$. Let us therefore recall some trivial but useful observations about $\lp{1}$-spaces.

If $X$ and $Y$ are two sets then it is well known that there is an isometric linear isomorphism $\lp{1}(X)\ptp\lp{1}(Y) \iso \lp{1}(X\times Y)$: see e.g. \cite{Ryan_TP}, Example 2.6. Moreover, if $E$ is a Banach space then every bounded linear map $T:\lp{1}(X) \to E$ restricts to a bounded function $X \to E$, where elements of $X$ are identified with the corresponding basis vectors in $\lp{1}(X)$. Conversely, any bounded function $X \to E$ extends uniquely to a bounded linear function $\lp{1}(X) \to E$. 

Therefore, if $X_1, \ldots, X_n$ are sets, every bounded function from $X_1\times\ldots\times X_n$ to $E$ extends uniquely to a bounded linear map from $\lp{1}(X_1)\ptp\ldots \ptp \lp{1}(X_n)$ to $E$.

These observations motivate the following definition, which is not standard but provides useful shorthand notation. Let $e_x$ denote the element of $\lp{1}(X)$ which takes the value 1 on $x$ and is 0 everywhere else.
\begin{defn}
Let $X_1, \ldots, X_n$ be sets. A \dt{primitive tensor} in $\lp{1}(X_1)\ptp \ldots\ptp\lp{1}(X_n)$ is an elementary tensor of the form
\[ e_{x_1} \tp \ldots \tp e_{x_n} \]
where $x_i \in X_i$ for $i=1, \ldots, n$. We shall often abuse notation by identifying a primitive tensor in $\lp{1}(X)\ptp\ldots\ptp\lp{1}(X)$ with an $n$-tuple of elements of $X$, and shall denote both by a bold letter: thus if $\vctR{x}=(x_1, \ldots, x_n) \in X^n$, we shall identify $\vctR{x}$ with $e_{x_1} \tp \ldots \tp e_{x_n}$.
\end{defn}
Since the linear span of $\{e_x \st x \in X\}$ is dense in $\lp{1}(X)$, the set of primitive tensors in $\lp{1}(X_1)\ptp \ldots\ptp\lp{1}(X_n)$ spans a dense linear subspace.

\begin{subsection}{Hochschild cohomology}
In this section we recall some basic facts about the (continuous) Hochschild cohomology of Banach algebras (see \cite {Hel_HBTA} or \cite{BEJ_CIBA} for more details). 

Let $A$ be a Banach algebra and $M$ a Banach $A$-bimodule. For each $n \geq 0$ let
\[ \begin{aligned}
\Ho{C}{n}(A,M) & \defeq M\ptp A^{\ptp n} \\
\Co{C}{n}(A,M) & \defeq \{ \text{bounded $n$-linear maps $A^n \to M$}\}
\end{aligned} \]
Elements of $\Ho{C}{n}(A,M)$ [resp. $\Co{C}{n}(A,M)$] are called \dt{continuous $n$-chains} [resp. \dt{$n$-cochains}] on $A$ 
with coefficients in $M$.

The Banach spaces $\Ho{C}{n}(A,M)$ and $\Co{C}{n}(A,M)$ fit into chain and cochain complexes, respectively, as follows:
\begin{subequations}
\begin{equation}\label{eq:hoch_ho}
0 \leftarrow \Ho{C}{0}(A,M) \lTo^{\bdy_0} \Ho{C}{1}(A,M) \lTo^{\bdy_1} \Ho{C}{2}(A,M) \lTo^{\bdy_2}  \ldots \end{equation}
\begin{equation}\label{eq:hoch_co}
 0 \rightarrow \Co{C}{0}(A,M) \rTo^{\dif_0} \Co{C}{1}(A,M) \rTo^{\dif_1} \Co{C}{2}(A,M) \rTo^{\dif_2} \ldots
\end{equation}
\end{subequations}
where the \dt{Hochschild boundary operator} $\bdy_n:\Ho{C}{n+1}(A,M) \to \Ho{C}{n}(A,M)$ is defined by
\begin{equation*}
\bdy_n(x\tp a_1\tp\ldots\tp a_{n+1}) =
 \left\{\begin{aligned}
& xa_1\tp a_2\tp\ldots\tp a_{n+1} \\
& + \sum_{j=1}^n (-1)^j x\tp a_1\tp\ldots\tp a_ja_{j+1}\tp\ldots\tp a_{n+1} \\
& + (-1)^{n+1} a_{n+1}x\tp a_1\tp\ldots\tp a_n
\end{aligned}\right.
\end{equation*}
and the \dt{Hochschild coboundary operator} $\dif_n:\Co{C}{n}(A,M) \to \Co{C}{n+1}(A,M)$ is defined by
\begin{equation*}
(\dif_n F)(a_1, \ldots, a_{n+1}) = \left\{\begin{aligned}
& a_1\cdot F(a_2, \ldots, a_{n+1}) \\
& + \sum_{j=1}^n (-1)^j F(a_1, \ldots, a_ja_{j+1}, \ldots, a_{n+1}) \\
& + (-1)^{n+1} F(a_1, \ldots, a_n)\cdot a_{n+1} \end{aligned}\right.
\end{equation*}
The $n$th homology group of the chain complex \eqref{eq:hoch_ho} is, by defintion, the quotient vector space
\[ \Ho{H}{n}(A,M) \defeq \frac{\Ker \bdy_{n-1}}{\Image \bdy_n} \]
and the $n$th cohomology group of the cochain complex \eqref{eq:hoch_co} is
\[ \Co{H}{n}(A,M) \defeq \frac{\Ker \dif_n}{\Image \dif_{n-1}} \]
We refer to the vector spaces $\Ho{H}{*}(A,M)$ and $\Co{H}{*}(A,M)$ as the \dt{Hochschild homology} and \dt{cohomology} groups, respectively, of $A$ with coefficients in $M$.

Of particular interest are the cases where the bimodule $M$ is either $A$ itself (with the bimodule action given by left and right multiplication in $A$) or $A'$ (with the canonically induced dual action). We refer to the spaces $\Ho{H}{n}(A,A)$ as the \dt{simplicial homology groups} of $A$ and to $\Co{H}{n}(A,A')$ as the \dt{simplicial cohomology groups} of $A$. We shall also use the notation $\SHo{H}{*}{A}$ and $\SCo{H}{*}{A}$ as abbreviations for the simplicial homology and cohomology groups respectively.

\end{subsection}

\end{section}
\begin{section}{Reduction to simplicial homology}
Under certain conditions, the simplicial homology of a commutative Banach algebra gives us information on Hochschild cohomology with symmetric coefficients. The precise statement requires us to consider the \dt{tensor product over $A$} of Banach $A$-modules: a full acount can be found in \cite{Hel_HBTA}, \S II.3 but we shall briefly summarise the required facts below.

If $X$ is a right Banach $A$-module and $Y$ a left Banach $A$-module, the Banach space $X\ptpR{A}Y$ is defined to be the quotient of $X\ptp Y$ by the closed linear span of the set
\[  \{xa\tp y - x\tp ay \st a\in A, x\in X, y \in Y \} \;; \]
and if $Z$ is a left Banach $A$-module, the Banach space $\lHom{A}(Y,Z)$ is defined to be the space of all bounded, $A$-module maps from $Y$ to $Z$, equipped with the usual operator norm.

It is important for the ensuing argument that both $X\ptpR{A}\blank$ and $\lHom{A}(\blank, Z)$ define additive functors from the category of left Banach $A$-modules to the category of Banach spaces: in particular, if
\[ \ldots \leftarrow C_{n-1} \leftarrow C_n \leftarrow \ldots \]
is a \emph{split exact} complex of left Banach $A$-modules and $A$-module maps, then the induced complexes $X\ptpR{A}C_*$ and $\lHom{A}(C_*, Z)$ of Banach spaces are both split exact.

\begin{lemma}\label{l:basechange}
Suppose $A$ is a commutative Banach algebra, and regard each Banach space $\Ho{C}{n}(A)=A\ptp A^{\ptp n}$ as a left Banach $A$-module via
\[ b\cdot(a_0\tp a_1\tp\ldots a_n) =ba_0\tp a_1\tp \ldots a_n \quad\quad(b,a_0, \ldots, a_n \in A). \]
Then for each $n \geq 0$ the Hochschild boundary map $\bdy_n:\Ho{C}{n+1}(A)\to\Ho{C}{n}(A)$ is a map of left Banach $A$-modules.

Moreover: if $A$ is also unital, and $M$ is a symmetric, unit-linked $A$-bimodule, there are isometric isomorphisms of chain complexes of Banach spaces
\[ \Ho{C}{*}(A,M) \iso M_R \ptpR{A} \SHo{C}{*}{A}  \]
and
\[ \Co{C}{*}(A,M) \iso {}_A{\rm Hom}( \SHo{C}{*}{A}, {}_LM) \]
where $M_R$ and ${}_LM$ are the right and left $A$-modules obtained from $M$ by restricting the action to one side.
\end{lemma}
\begin{proof}
The proof that $\bdy_n:\Ho{C}{n+1}(A)\to\Ho{C}{n}(A)$ is an $A$-module map is a routine calculation (it is crucial here that $A$ is commutative). The statement about isomorphism of chain complexes is a simple extension of the following observation: when $A$ is unital and $M$ is symmetric, then for any Banach space $E$ there are isometric isomorphisms of Banach spaces
\[ M\ptp E \iso M \ptpR{A} (A\ptp E)\quad\text{ and }\quad {\mc L}(E,M)\iso {}_A{\rm Hom} (A\ptp E, {}_LM) \]
\end{proof}

Next we state the main result of this section, for which the author knows no exact reference in the literature, but which is probably not a new observation.
\begin{propn}\label{p:UCT}
Let $A$ be a unital, commutative Banach algebra that is isomorphic as a Banach space to $\lp{1}(\Omega)$ for some indexing set $\Omega$. The following are equivalent:
\begin{itemize}
\item[(a)] $\SCo{H}{n}{A}=0$ for all $n \geq 1$;
\item[(b)] $\SHo{H}{n}{A}=0$ for all $n \geq 1$;
\item[(c)] the complex of left Banach $A$-modules
\[0 \leftarrow \SHo{C}{1}{A} \lTo^{\bdy_1} \SHo{C}{2}{A} \lTo^{\bdy_2} \ldots \]
splits in the category of Banach spaces and bounded linear maps;
\item[(d)] the complex of left $A$-modules
\[0 \leftarrow \SHo{C}{1}{A} \lTo^{\bdy_1} \SHo{C}{2}{A} \lTo^{\bdy_2} \ldots \]
splits in the category of left Banach $A$-modules;
\item[(e1)] $\Co{H}{n}(A,M)=0$ for all $n \geq 1$ and any \emph{symmetric, unit-linked} $A$-bimodule $M$;
\item[(e2)] $\Ho{H}{n}(A,M)=0$ for all $n \geq 1$ and any \emph{symmetric, unit-linked} $A$-bimodule $M$.
\end{itemize}
\end{propn}

\begin{proof}
Clearly $(e1)\implies(a)$ and $(e2)\implies (b)$. We shall show that
\[ (a)\iff(b)\implies(c)\implies(d)\implies (e1) \& (e2) \]
which will prove the theorem.

\vspace{0.5em}\noindent\underline{\it Proof of $(a)\iff (b)$}: 
The Hochschild cochain complex
\[  0 \to \SCo{C}{0}{A} \rTo^{\dif_0} \SCo{C}{1}{A} \rTo^{\dif_1}  \ldots \]
is the dual of the Hochschild chain complex
\[ 0 \leftarrow \SHo{C}{0}{A} \lTo^{\bdy_0} \SHo{C}{1}{A} \lTo^{\bdy_1}  \ldots \]
that is, $\bdy_n'=\dif_n$ for all $n$. Since $A$ is commutative, $\bdy_0=0$ and $\dif_0=0$. Hence (a) is equivalent to exactness of the complex
\begin{subequations}
\begin{equation}\label{eq:FRED1}
0 \to  \SCo{C}{1}{A} \rTo^{\dif_1} \SCo{C}{2}{A} \rTo^{\dif_2} \ldots
\end{equation}
while (b) is equivalent to exactness of the complex
\begin{equation}\label{eq:FRED2}
 0 \leftarrow \SHo{C}{1}{A} \lTo^{\bdy_1} \SHo{C}{2}{A} \lTo^{\bdy_2} \ldots
\end{equation}
\end{subequations}
Since \eqref{eq:FRED2} is the dual of \eqref{eq:FRED1}, exactness of one of the complexes implies exactness of the other by astandard duality argument (see \cite{BEJ_CIBA} Corollary 1.3, or \cite{Hel_HBTA} Proposition II.5.29).

\vspace{0.5em}\noindent\underline{\it Proof of $(b)\implies(c)$}: 
We have seen that $(b)$ holds if and only if the chain complex
\[0 \leftarrow \SHo{C}{1}{A} \lTo^{\bdy_1} \SHo{C}{2}{A} \lTo^{\bdy_2} \ldots \]
is exact. Recall that $\SHo{C}{n}{A}\iso A\ptp A^{\ptp n} \iso \lp{1}(\Om^{n+1})$ for each $n$. Since $\bdy_1$ is surjective, the well-known lifting property of $\SHo{C}{1}{A}=\lp{1}(\Om^2)$ with respect to open mappings allows us to find a bounded linear map $s_1:\SHo{C}{1}{A}\to\SHo{C}{2}{A}$ such that $\bdy_1s_1={\sf id}$. Then ${\sf id}- s_1\bdy_1$ is a bounded linear projection of $\SHo{C}{2}{A}$ onto $\SHo{Z}{2}{A}=\Ker(\bdy_1)=\Image(\bdy_2)$.

Using the lifting property of $\SHo{C}{2}{A}=\lp{1}(\Om^3)$ with respect to the surjection $\SHo{C}{3}{A}\rTo^{\bdy_2}\SHo{Z}{2}{A}$, we can lift the map ${\sf id}-s_1\bdy_1$ to a bounded linear map $s_2:\SHo{C}{2}{A}\to\SHo{C}{3}{A}$: by construction, $\bdy_2s_2={\sf id}-s_1\bdy_1$.

Continuing in this way, we inductively construct bounded linear maps $s_n:\SHo{C}{n}{A}\to\SHo{C}{n+1}{A}$ for each $n \geq 1$, such that $s_n\bdy_n+\bdy_{n+1}s_{n+1}={\sf id}$ for all $n$.

\vspace{0.5em}\noindent\underline{\it Proof of $(c)\implies(d)$}: 
Condition $(c)$ states that the complex of Banach $A$-modules
\[ 0 \leftarrow \SHo{C}{1}{A} \lTo^{\bdy_1} \SHo{C}{2}{A} \lTo^{\bdy_2} \ldots \]
splits in the category of Banach spaces and bounded linear maps. But for each $n \geq 0$, $\SHo{C}{n}{A}\iso A\ptp A^{\ptp n}$ is projective (\cite{Hel_HBTA}, Defn III.1.13) as a left Banach $A$-module, since $A$ is \emph{unital}. Therefore, if $(c)$ holds, the projectivity of each $\SHo{C}{n}{A}$ may be used to inductively construct continuous $A$-module maps  $\sig_n:\SHo{C}{n}{A}\to\SHo{C}{n+1}{A}$, $n \geq 1$, such that $\sig_n\bdy_n+\bdy_{n+1}\sig_{n+1}={\sf id}$ for all $n$. The argument is essentially the same as that used to prove that $(b)\implies(c)$.

\vspace{0.5em}\noindent\underline{\it Proof of $(d)\implies(e1)\&(e2)$}: 
By Lemma \ref{l:basechange}, for each $n$ we have
\[ \begin{aligned}
& \Ho{H}{n}(A,M) & = & H_n\left( M_R \ptpR{A} \SHo{C}{*}{A} \right) \\
\text{ and } & \Co{H}{n}(A,M) & = & H^n\left( {}_A{\rm Hom}( \SHo{C}{*}{A}, {}_LM)\; \right) \end{aligned} \]
where $M_R$ and ${}_LM$ are the right and left $A$-modules obtained from $M$ by restricting the action to one side. Now if $(d)$ holds then the Banach space complexes $M_R \ptpR{A} \SHo{C}{*}{A}$ and ${}_A{\rm Hom}( \SHo{C}{*}{A}, {}_LM )$ are both (split) exact in degrees $\geq 1$, and so $(e1)$ and $(e2)$ must hold.
\end{proof}

Before going on to apply this result, we make some remarks to put it in context: these are not required for the main result of this paper, but may serve as motivation for Proposition \ref{p:UCT}.

\begin{subequations}
Specifically, Proposition \ref{p:UCT} is modelled on a special case of the following result for ``purely algebraic'' Hochschild cohomology (see \cite{Weibel}, Ch. 9 for the necessary definitions).
\begin{propn}
Let $k$ be a commutative ring, and let $R$ be a commutative unital $k$-algebra \emph{which is $k$-projective}. Then for any unit-linked, symmetric $R$-bimodule $X$, there is a spectral sequence
\begin{equation}\label{eq:base-change-ss}
\Ext_R^p \left( H_q(R,R),\,X\right) \underset{p}{\Longrightarrow} H^{p+q}(R,X)
\end{equation}
\end{propn}
Although we cannot find an exact reference for \eqref{eq:base-change-ss}, it follows easily from the isomorphisms $H_q(R,R)\iso \Tor^{R^e}_q (R,R)$, $H^n(R,X) \iso \Ext_{R^e}(R,X)$ (\cite{Weibel}, Corollary 9.1.5) and a change-of-rings spectral sequence for $\Ext$ (see, for example, \cite{Rot_HA} Theorem 11.6.5).

 In particular, if $H_q(R,R)=0$ for all $q\geq 1$, the spectral sequence collapses to give
\begin{equation}\label{eq:alg-simp-triv}
 H_n(R,X) \iso \Ext_R^n\left( H_0(R,R),\,X\right) = 0 \quad\forall\,n\geq 1 \end{equation}
where the last equality holds because $H_0(R,R)=R$ is $R$-projective.
\end{subequations}

The point is that Equation \eqref{eq:alg-simp-triv} is the algebraic analogue of the implication $(b)\implies(e1)$ of Proposition \ref{p:UCT}. We hope that the use of Proposition \ref{p:UCT} in this paper may lead to further work in finding effective, Banach algebraic analogues of the change-of-rings spectral sequences.
\end{section}

\begin{section}[Statement of main theorem]{Statement of main theorem and initial reductions}\label{s:state_thm}
\begin{defn}
A Banach algebra $A$ is \dt{simplicially trivial} if $\SCo{H}{n}{A}=0$ for all $n \geq 1$.
\end{defn}
Since a Banach algebra $A$ is amenable if and only if $\Co{H}{n}(A,X')=0$ for all $A$-bimodules $X$ and all $n \geq 1$ (\cite{Bons-Dunc} \S44 Proposition 6), every amenable Banach algebra is simplicially trivial. However, simplicial triviality is in general much weaker than amenability. It was observed at the start that if $S$ is an infinite semilattice then by \cite{Dunc-Nam}, Theorem 10, $\A{S}$ cannot be amenable. In contrast we present the following theorem, which is the main result of this paper.

\begin{thm}\label{t:result}
Let $S$ be a unital semilattice. Then $\A{S}$ is simplicially trivial.
\end{thm}

\begin{coroll}\label{c:symm_Hoch}
Let $S$ be a semilattice. Then the Hochschild cohomology groups $\Co{H}{n}(\A{S},M)$ are zero for all $n\geq 1$ and all symmetric $\A{S}$-bimodules $M$.
\end{coroll}
\begin{proof}[Proof that Theorem \ref{t:result}$\implies$ Corollary \ref{c:symm_Hoch}]
We use Proposition \ref{p:UCT}, together with some cohomological machinery to get round the problem that $\A{S}$ might not be unital.

If $A$ is any Banach algebra, with forced unitisation denoted by $\fu{A}$, then each Banach $A$-bimodule $X$ may be regarded as a unit-linked, Banach $\fu{A}$-bimodule $X_1$, where the action of $\fu{A}$ on $X_1$ is given by
\[ (a+\lm{\sf 1})\cdot x \defeq ax+\lm x \quad;\quad x\cdot(a+\lm{\sf 1}) \defeq xa+\lm a \quad\quad(a,\in A,\lm\in\Cplx).\]
There is then a canonical restriction map
\[ \Co{C}{n}(\fu{A}, X_1) \to \Co{C}{n}(A,X) \]
which induces an isomorphism of cohomology groups for each $n$ (\cite{Hel_HBTA}, Exercise III.4.10; see also \cite{BEJ_CIBA}, Ch. 1).

In particular, $\Co{H}{n}(\fu{\A{S}},M_1) \iso \Co{H}{n}(\A{S},M)$ for all $n$. Now since $\fu{(\A{S})}=\A{(\fu{S})}$,  Theorem \ref{t:result} implies that $\fu{\A{S}}$ is simplicially trivial. Invoking Proposition \ref{p:UCT} we deduce that $\Co{H}{n}(\fu{\A{S}},M_1)=0$ for all $n\geq 1$, and this completes the proof.
\end{proof}

Theorem \ref{t:result} will be deduced from a stronger statement (Theorem \ref{t:mainthm} below). The proof requires some preliminaries, which are given in the next two sections.
\end{section}

\begin{section}[The unital, finite free case]{A contracting homotopy for the unital, finite free case}\label{s:finite_free}
In this section, fix a finite set $X$ and let $F=2^X$ be the free unital semilattice generated by $X$, as defined in Example \ref{eg:finfree}. For notational convenience we write $e_i$ rather than $e_{\{i\}}$ for the element of $\lp{1}(F)$ that takes the value 1 on the element $\{i\}$ and is 0 everwhere else: note that by the way multiplication in $F$ and $\A{F}$ is defined,
\[ e_J = \prod_{i \in J} e_i \quad\quad\text{ for all $J \subseteq X$.} \]

\begin{lemma}\label{l:fin_free_split}
There exist linear maps $s_n:\Ho{C}{n}(\A{F}) \to \Ho{C}{n+1}(\A{F})$, $n \geq 0$, such that ${\bdy_ns_n + s_{n-1}\bdy_{n-1}={\sf id}}$ for all $n \geq 1$.
\end{lemma}

We emphasise that for our application to the proof of Theorem \ref{t:result}, explicit formulas for the maps $s_n$ are not required. However finding explicit $s_n$ is straightforward and so we give a sharper form of Lemma \ref{l:fin_free_split} below.
To do this we introduce distinguished elements of the algebra $\A{F}$: for every $J \subseteq F$ let
\[ u_J = \prod_{i \in J} e_i \prod_{k \in X \setminus J} (e_{\emptyset}-e_k) \;\in\A{F}\;. \]

\begin{lemma}\label{l:fin_free_explicit}
For $n \geq 1$ define {$s_n: \Ho{C}{n}(\A{F}) \to \Ho{C}{n+1}(\A{F})$} by
\begin{equation}
s_n (a_0\tp a_1\tp\ldots\tp a_n) \defeq \sum_{J \subseteq X} a_0u_J\tp u_J \tp a_1\tp \ldots\tp a_n\quad\quad{(a_1, \ldots, a_n \in \A{F})}
\end{equation}
and let $s_0=0$. Then $\bdy_n s_n + s_{n-1}\bdy_{n-1}={\sf id}$ for all $n \geq 1$, and $\norm{s_n}\leq 5^{\abs{X}}$ for all $n$.
\end{lemma}

The required properties of $u_J$ are collected in the following lemma.
\begin{lemma}\label{l:diagonal}
Let $J \subseteq X$. Then:
\begin{itemize}
\item[(i)] $\norm{u_J}\leq 2^{\abs{X}-\abs{J}}$;
\item[(ii)] $au_J\tp u_J = u_J\tp u_Ja$ for all $a \in \A{F}$.
\item[(iii)] $\sum_{J \subseteq X} u_J^2 =e_{\emptyset}$;
\end{itemize}
\end{lemma}
\begin{proof}
The norm estimate $(i)$ is trivial.

Since $F$ is commutative, for each $i\in X$ a direct computation yields
\[  e_i u_J = \left\{\begin{aligned} 0 & \text{ if $i \notin J$} \\ u_J  & \text{ if $i \in J$} \end{aligned}\right\} = u_J e_i \;.\]
It follows immediately that $e_iu_J\tp u_J = u_J\tp u_Je_i$ for all $i \in X$, and since the $e_i$ span $\A{F}$ as a vector space we deduce by linearity that $au_J\tp u_J=u_J\tp u_Ja$. Thus (ii) is proved.

Finally: since $F$ is commutative $u_J$ is the product of commuting idempotents and is therefore itself an idempotent. Hence
\[ \sum_{J \subseteq X} u_J^2 = \sum_{J \subseteq X} u_J = \prod_{i \in X} \bigl(e_i + (e_{\emptyset}-e_i) \bigr) = e_{\emptyset} \;.\]
and (iii) is proved.
\end{proof}

\begin{rem}
Those familiar with Hochschild cohomology will recognise that we are exploiting the existence of a ``diagonal for the algebra $\A{F}$''. To put this in context: whenever $A$ is an algebra with a diagonal $\Delta$ and $M$ is an $A$-bimodule, standard techniques from homological algebra can be applied to construct from $\Delta$ an explicit contracting homotopy of the chain complex $\Ho{C}{*}(A,M)$. However, in order to appeal to known results one must first prove that $\A{F}$ has a diagonal, and such a proof is not significantly shorter than the direct approach taken here.
\end{rem}

\begin{proof}[Proof of Lemma \ref{l:fin_free_explicit}]
The proof that 
\[ \bdy_n s_n + s_{n-1}\bdy_{n-1}={\sf id} \]
is a direct computation using properties (ii) and (iii) in Lemma \ref{l:diagonal}. The estimate on $\norm{s_n}$ follows from the inequalities
\[ \norm{s_n (a_0\tp a_1\tp\ldots\tp a_n)}\leq \sum_{J \subseteq X} \norm{u_J}^2 \norm{a_0}\ldots\norm{a_n} \]
and
\[ \sum_{J\subseteq X} \norm{u_J}^2 \leq \sum_{J \subseteq X} 4^{\abs{X}-\abs{J}} = \sum_{k=0}^{\abs{X}} {\abs{X} \choose k} 4^{\abs{X}-k} =  5^{\abs{X}} \]
where the last estimate uses the bound established in Lemma \ref{l:diagonal} (i).
\end{proof}

For fixed $n$, the splitting map $s_n$ that is given in Lemma \ref{l:fin_free_explicit} clearly depends on the underlying set $X$ which generates $F$, and there is no reason to believe that there is a uniform bound on the norm of $s_n$ as the size of $X$ increases. (In fact, there can be no such bound: if there were, a straightforward exhaustion argument would allow one to construct a bounded approximate diagonal for $\A{2^{\Nat}}$, giving a contradiction since by \cite{Dunc-Nam} $\A{2^{\Nat}}$ is not amenable.) Nevertheless, in Section \ref{s:mainthm} we shall use the maps $s_n$ to inductively construct \emph{natural} splitting maps, and it is these natural maps which split the simplicial chain complex of any $\lp{1}$ semilattice algebra.

Here the word ``natural" must be made precise, in terms of functors between categories. This will be done in the next section.
\end{section}

\begin{section}{Natural transformations and formal substitution}\label{s:formalsub}
We shall only need the basic language of category theory, for which a standard reference is \cite{Mac_CWM}.

The categories $\SGp$ and $\BAlg_1$ have already been introduced. We shall also need to consider the following:
\begin{itemize}
\item $\SGp_*$ -- the category of monoids and monoid homomorphisms;
\item $\SLat_*$ -- the full subcategory of $\SGp_*$ whose objects are unital semilattices;
\item $\BAlg$ -- the category of Banach algebras and continuous algebra homomorphisms;
\item $\Ban$ -- the category of Banach spaces and continuous linear maps.
\end{itemize}

\begin{lemma}\label{l:hoch_is_nat}
For any $n \geq 0$, $\SHo{C}{n}{\blank}$ defines a functor from $\BAlg$ to $\Ban$. Moreover, the Hochschild boundary map $\bdy_n : \SHo{C}{n+1}{\blank} \to \SHo{C}{n}{\blank}$ defines a natural transformation between functors.
\end{lemma}

\begin{proof}
The first part is immediate once one recalls that the underlying Banach space of $\SHo{C}{n}{A}$ is $A^{\ptp n+1}$, and that for each $m$ taking the $m$th projective tensor power is a functorial operation on Banach spaces. The second part is an easy calculation, which depends on the fact that each face map
\[ \partial_j: \SHo{C}{n+1}{A}\to\SHo{C}{n}{A}\quad;\quad a_0 \tp \ldots \tp a_n \mapsto a_0\tp\ldots\tp a_ja_{j+1}\tp\ldots\tp a_n \]
is natural in $A\in\BAlg$.
\end{proof}

If $A \rTo^f B$ in $\BAlg_1$ we shall write $f^{\tp n+1}$ rather than $\SHo{C}{n}{f}$. We shall abuse notation further by writing $\bdy^H_n$ for $\bdy^{\A{H}}_n$ whenever $H$ is a semilattice. Then for each $n$, $\SHo{C}{n}{\A{\blob}}$ may be considered as a functor from $\SLat_*$ to $\Ban$ -- more precisely, as the composite of the functors
\[ \SLat_* \rTo \SGp_* \rTo^{\A{\blob}} \BAlg \rTo^{\SHo{C}{n}{\blank}} \Ban  \]

\begin{inline}{Important convention.} Henceforth, whenever we refer to one of the functors $\SHo{C}{n}{\A{\blob}}$ we shall always mean that it is a functor from $\SLat_*$ to $\Ban$.

This allows us to regard the family $(\bdy_n^H)_{H \in \SLat_*}$ as a natural transformation from $\SHo{C}{n+1}{\A{\blob}}$ to $\SHo{C}{n}{\A{\blob}}$
\begin{diagram}
\SLat & \rTo^{\SHo{C}{n+1}{\A{\blob}}} & \Ban \\
\dEq & \dDbl_{\bdy_n} & \dEq \\
\SLat & \rTo_{\SHo{C}{n}{\A{\blob}}} & \Ban
\end{diagram}
\end{inline}

\subsection*{Formal substitution}
Let $j \geq 1$ and let $F$ be the free unital semilattice on $j+1$ generators $f_0, \ldots, f_j$. Given a unital semilattice $S$ and $\vctR{x}=(x_0, \ldots, x_j) \in S^{j+1}$, let
\[ \subst{x}(f_i) \defeq x_i \quad\quad(i=0,1,\ldots, j) \]
Because each $x_i$ is an idempotent and the $f_i$ are free generators for $F$, this uniquely defines a homomorphism of monoids from $F$ to $S$, and thus defines a (unique) morphism {$\A{F} \rTo^{\subst{x}} \A{S}$} in $\BAlg_1$ with the property that $\subst{x}(e_{f_i})=e_{x_i}$ for $i=0,1, \ldots j$. This in turn yields the following ``substitution property''
\begin{equation}\label{eq:sub-prop}
 \subst{x}^{\tp j+1}(\vctR{f}) = \subst{x}(e_{f_0})\tp\ldots\subst{x}(e_{f_n})=\vctR{x}
\end{equation}
which will be crucial to our main construction. (Informally, $\subst{x}$ stands for ``formal substitution of $x_i$ for each occurence of $f_i$''.)

Applying Lemma \ref{l:hoch_is_nat} gives the following commutative diagram in $\Ban$:
\begin{diagram}
\ldots &\lTo & \SHo{C}{n-1}{\A{F}} & \lTo^{\bdy^F_{n-1}} & \SHo{C}{n}{\A{F}} &  \lTo^{\bdy^F_n} & \SHo{C}{n+1}{\A{F}} & \lTo \ldots\\
& & \dTo^{\subst{x}^{\tp n}} & & \dTo_{\subst{x}^{\tp n+1}} & & \dTo_{\subst{x}^{\tp n+2}} & \\
\ldots & \lTo & \SHo{C}{n-1}{\A{S}} & \lTo_{\bdy^S_{n-1}} & \SHo{C}{n}{\A{S}} & \lTo_{\bdy^S_n} & \SHo{C}{n+1}{\A{S}} & \lTo \ldots
\end{diagram}

\end{section}

\begin{section}{The main technical theorem}\label{s:mainthm}
In this section we shall use the preceding ideas to prove Theorem \ref{t:result}. Our proof technique turns out to yield something stronger:
\begin{thm}\label{t:mainthm}
Let $H$ be a unital semilattice and $\A{H}$ its $\lp{1}$-convolution algebra.  Then for each $n \in \Zahl_+$ there exist bounded linear maps $\sig^H_n: \SHo{C}{n}{\A{H}} \to \SHo{C}{n+1}{\A{H}}$, natural in $H\in\SLat_*$, such that
\[ \bdy^H_n\sig^H_n+\sig^H_{n-1}\bdy^H_{n-1}={\sf id} \quad\text{for all $n\geq 1$}\]
In particular, the simplicial Hochschild chain complex $\SHo{C}{n}{\A{H}}$ is split exact in degrees 1 and above.
\end{thm}
The proof of this stronger result (Theorem \ref{t:mainthm}) will take up the rest of this section.
The idea is to construct the natural transformations $\sig_j: \SHo{C}{j}{\A{\blob}} \to \SHo{C}{j+1}{\A{\blob}}$ \emph{recursively}: the naturality assumption is needed in the inductive hypothesis for our construction to work.

We isolate the inductive step as a lemma in its own right.
\begin{lemma}\label{l:inductive_step}
Let $j \geq 1$. Suppose that there exists a natural transformation ${\sig_{j-1}: \SHo{C}{j-1}{\A{\blob}} \to \SHo{C}{j}{\A{\blob}} }$ such that
\begin{equation}\label{eq:inductive_hyp}
 \bdy^S_{j-1}\sig^S_{j-1}\bdy^S_{j-1}=\bdy^S_{j-1}
 \quad\quad(S\in\SLat_*)
\end{equation}
Then there exists a natural transformation $\sig_j: \SHo{C}{j}{\A{\blob}}\to \SHo{C}{j+1}{\A{\blob}}$ such that ${\bdy^S_j\sig^S_j+ \sig^S_{j-1}\bdy^S_{j-1}} = {\sf id}$ for every $S\in\SLat_*$.
\end{lemma}

\begin{proof}
Let $F=2^{[j+1]}$ be the free unital semilattice generated by $j+1$ idempotents $f_0, \ldots, f_j$. We denote the $(j+1)$-tuple $(f_0,f_1, \ldots, f_j)$ and the associated primitive tensor in $\A{F}^{\tp j+1}$ by $\vctR{f}$. By Lemma \ref{l:fin_free_split}, there exist bounded linear maps $s^j_n: \SHo{C}{n}{\A{F}} \to \SHo{C}{n+1}{\A{F}}$ such that
\[ s^j_{n-1}\bdy^F_{n-1}+\bdy^F_ns^j_n ={\sf id} \]
for all $n \in \Nat$.

In particular, $s^j_{j-1}\bdy^F_{j-1}+\bdy^F_js^j_j ={\sf id}$, and so
\begin{equation}\label{eq:bad_split}
\begin{split}
{\sf id}-\sig^F_{j-1}\bdy^F_{j-1}
& = (\bdy^F_js^j_j+s^j_{j-1}\bdy^F_{j-1})({\sf id}-\sig^F_{j-1}\bdy^F_{j-1}) \\
& = \bdy^F_j s^j_j ({\sf id}-\sig^F_{j-1}\bdy^F_{j-1})
 + s^j_{j-1}(\bdy^F_{j-1}-\bdy^F_{j-1}\sig^F_{j-1}\bdy^F_{j-1}) \\
& =  \bdy^F_j s^j_j ({\sf id}-\sig^F_{j-1}\bdy^F_{j-1})
\end{split}
\end{equation}
where the last equality follows from the starting hypothesis \eqref{eq:inductive_hyp}. Applying both sides of \eqref{eq:bad_split} to the primitive tensor $\vctR{f}$ yields the following equation in $\A{F}^{\tp j+1}$:
\begin{equation}\label{eq:generic}
 \vctR{f}-\sig^j_{j-1}\bdy^F_{j-1}\vctR{f}
=  \bdy^F_j s^j_j \left( \vctR{f} -\sig^F_{j-1}\bdy^F_{j-1}\vctR{f} \right)
\end{equation}

Now define ${\sf w} \in \SHo{C}{j+1}{\A{F}}$ by
\[\higheq{ {\sf w} \defeq s^j_j (\vctR{f}-\sig^F_{j-1}\bdy^F_{j-1}\vctR{f}) } \]
and observe that by construction
\begin{equation}\label{eq:formal_identity}
\higheq{ \bdy^F_j({\sf w})\;=\; \vctR{f} - \sig^F_{j-1}\bdy^F_{j-1}(\vctR{f}) }
\end{equation}
Intuitively, Equation \eqref{eq:formal_identity} is a ``formal identity'' and so will hold if $f_0,f_1, \ldots, f_j$ are replaced with arbitrary elements of an arbitrary unital semilattice. To make this argument precise, we use the substitution morphisms $\subst{x}$ which were introduced in Section \ref{s:formalsub}.

Let $S\in\SLat_*$ be a unital semilattice, and recall that ${\sf w}$ does not depend on $S$ (it does however depend on the map $s^j_j$, and on the presumed existence of the maps $\sig^F_{j-1}$). Define $\sig^S_j: \SHo{C}{j}{\A{S}} \to \SHo{C}{j+1}{\A{S}}$ on \emph{primitive} tensors by
\[ \sig^S_j(\vctR{x}) \defeq \subst{x}^{\tp j+2} ({\sf w}) \quad\quad(x_0, \ldots, x_j \in S) \]
and extend by linearity: since $\subst{x}$ is a contractive linear map, so is $\subst{x}^{\tp j+2}$ and therefore
\[  \norm{\sig^S_j} = \sup_{\vctR{x} \in S^{j+1}} \norm{\sig^S_j(\vctR{x})} = \sup_{\vctR{x} \in S^{j+1}}  \norm{\subst{x}^{\tp j+2} ({\sf w})} \leq\norm{\sf w}_1 \]
Thus for each semilattice $S$, $\sig^S_j$ is a well-defined, bounded linear map (with norm less than or equal to $\norm{\sf w}_1$).

It remains to show that the family $(\sig^S_j)_{S \in \SLat_*}$ has the desired properties. First we show that $\bdy^S_j\sig^S_j + \sig^S_{j-1}\bdy^S_{j-1}={\sf id}$. By linearity and continuity it suffices to show that
\[ (\bdy^S_j\sig^S_j + \sig^S_{j-1}\bdy^S_{j-1})(\vctR{x})=\vctR{x} \]
for any primitive tensor $\vctR{x}$ in $\A{S}^{\ptp j+1}$. This is done using the following calculation, which is essentially a chase round the following commutative diagram:
\begin{diagram}
\SHo{C}{j-1}{\A{F}} & & \LR{\bdy^F_{j-1}}{\sig^F_{j-1}} & & \SHo{C}{j}{\A{F}} & &  \lTo^{\bdy^F_j} & & \SHo{C}{j+1}{\A{F}} \\
\dTo^{\subst{x}^{\tp j}} & & & & \dTo_{\subst{x}^{\tp j+1}} & & & & \dTo_{\subst{x}^{\tp j+2}} \\
\SHo{C}{j-1}{\A{S}} & & \LR{\bdy^S_{j-1}}{\sig^S_{j-1}} & & \SHo{C}{j}{\A{S}} & &\lTo^{\bdy^S_j} & & \SHo{C}{j+1}{\A{S}}
\end{diagram}
Fix $(x_0, \ldots, x_j) \in S^{j+1}$ and denote the corresponding primitive tensor in $\SHo{C}{j}{\A{F}}$ by $\vctR{x}$. Then
\[ \begin{aligned}
\bdy^S_j\sig^S_j(\vctR{x}) & = \bdy^S_j \subst{x}^{\tp j+2} ({\sf w})  &  \quad\text{ (definition of $\sig^S_j(\vctR{x})$)} \\
& = \subst{x}^{\tp j+1}\bdy^F_j({\sf w}) & \quad\text{ (naturality of $\bdy_j$)} \\
& = \subst{x}^{\tp j+1}\left( \vctR{f} - \sig^F_{j-1}\bdy^F_{j-1}(\vctR{f})\right)  & \quad\text{ (Equation \eqref{eq:formal_identity})} \\
& = \subst{x}^{\tp j+1}(\vctR{f}) -  \subst{x}^{\tp j+1}\sig^F_{j-1}\bdy^F_{j-1} (\vctR{f}) & \quad\text{ (linearity of $\subst{x}^{\tp j+1}$)} \\
& = \subst{x}^{\tp j+1}(\vctR{f}) - \sig^S_{j-1}\subst{x}^{\tp j}\bdy^F_{j-1}(\vctR{f}) & \quad\text{ (naturality of $\sig_{j-1}$)} \\
& = \subst{x}^{\tp j+1}(\vctR{f}) - \sig^S_{j-1}\bdy^S_{j-1}\subst{x}^{\tp j+1}(\vctR{f}) & \quad\text{ (naturality of $\bdy_{j-1}$)} \\
& = \vctR{x} - \sig^S_{j-1}\bdy^S_{j-1}(\vctR{x}) &  \quad\text{ (substitution property \eqref{eq:sub-prop} of $\subst{x}$)}
\end{aligned} \]
as required.

Finally, we need to show that $\sig_j$ is a natural transformation, i.e. that whenever ${\tht: H \to K}$ is a morphism in $\SLat_*$ we have a commuting square
\begin{diagram}
\SHo{C}{j}{\A{H}} & & \rTo^{\sig^H_j} & & \SHo{C}{j+1}{\A{H}} \\
\dTo^{\tht^{\tp j+1}} & & & & \dTo_{\tht^{\tp j+2}} \\
\SHo{C}{j}{\A{K}} & & \rTo_{\sig^K_j} & & \SHo{C}{j+1}{\A{K}}
\end{diagram}
By linearity and continuity, it suffices to check this on primitive tensors. Let $\vctR{a}$ be the primitive tensor corresponding to the $(j+1)$-tuple $(a_0,\ldots,a_j)$, where each $a_i \in H$. Writing $b_i$ for $\tht(a_i)\in K$, we have
\[ \sig^K_j\tht^{\tp j+1} (\vctR{a}) = \sig^K_j(\vctR{b}) = \subst{b}^{\tp j+2}({\sf w}) \]
while
\[ \tht^{\tp j+2}\sig^H_j(\vctR{a}) = \tht^{\tp j+2}\subst{a}^{\tp j+2}({\sf w}) = \left( \tht\subst{a} \right)^{\tp j+2}({\sf w}) \]
Since $\tht\subst{a} = \subst{b}$,
\[ \tht^{\tp j+2}\sig^H_j(\vctR{a}) = \subst{b}^{\tp j+2}({\sf w}) = \sig^K_j\tht^{\tp j+1} (\vctR{a}) \]
as required.
\end{proof}

\begin{proof}[Proof of Theorem \ref{t:mainthm}]
We shall prove the following:

\vspace{0.5em}\noindent{\bf Claim.} There exist natural transformations $\sig_m:\SHo{C}{m}{\A{\blob}} \to \SHo{C}{m+1}{\A{\blob}}$, $m \geq 0$, such that
\[  \sig^S_{n-1}\bdy^S_{n-1}+\bdy^S_n\sig^S_n={\sf id} \quad\quad\text{for all $n \geq 1$} \]
\vspace{0.5em}

The proof is by induction on $m$.

Observe that when $S\in\SLat_*$ $\A{S}$ is a commutative algebra, which implies that $\bdy_0^{S}:\SHo{C}{1}{\A{S}}\to\SHo{C}{0}{\A{S}}$ is the zero map. Therefore, if one takes $\sig_0=0$ then the conditions of Lemma \ref{l:inductive_step} are satisfied for $j=1$, and so by the lemma there exists a natural transformation $\sig_1$ such that
\[ \sig^S_0\bdy^S_0+\bdy^S_1\sig^S_1 ={\sf id} \;. \]

Suppose now that for some $j \geq 2$, there exist natural transformations\\ $\sig_0, \ldots, \sig_{j-1}$ such that 
\[ \sig^S_{i-1}\bdy^S_{i-1}+\bdy^S_i\sig^S_i ={\sf id} \]
whenever $1\leq i \leq j-1$ and $S\in\SLat_*$. Then in particular, $\bdy^S_{j-1}\sig^S_{j-1}\bdy^S_{j-1} = \bdy^S_{j-1}$. Applying Lemma \ref{l:inductive_step} we deduce that there exists a natural transformation $\sig_j$ such that
\[ \sig^T_{j-1} \bdy^T_{j-1}+\bdy^T_j\sig^T_j={\sf id} \quad\quad\text{ ($T\in\SLat_*$)}\]
This completes the inductive step; hence the claim, and with it Theorem \ref{t:mainthm}, is proved.
\end{proof}

\end{section}

\begin{section}{Closing thoughts}
It is natural to wonder if the proof technique set out in this article can be extended to other examples of commutative semigroups. Unfortunately, the proof relies on being able to find \emph{natural} splitting maps in some starting degree; without such a starting point the inductive argument of \S\ref{s:mainthm} is of no use. (In proving Lemma \ref{l:inductive_step}, it is the naturality assumption that allows us to transfer our understanding of the finite \emph{free} case to the general case.)

It is possible to obtain partial generalisations to algebras of the form $\A{S}$ where $S$ is a {\it Clifford semigroup}: this is part of work in progress (\cite{YC_Clifford-pre}).
\end{section}

\section*{Acknowledgements}
This work was carried out during the author's PhD studies at the University of Newcastle, which were supported by an EPSRC grant. The author would like to thank his supervisor Michael C.\ White for suggesting this line of enquiry and providing much support during its pursuit. He also acknowledges helpful feedback from Matthew Daws and Zinaida Lykova on early versions of this material, and the suggestions for improvement made by the referee on a first draft of this article.


\end{document}